\documentclass[letterpaper, 10 pt, conference]{ieeeconf}  

\IEEEoverridecommandlockouts                              

\usepackage{graphics} 
\usepackage{epsfig} 
\usepackage{amsmath} 
\usepackage{amssymb}  
\usepackage{cite}
\usepackage{xcolor}
\usepackage[lined,ruled,linesnumbered,commentsnumbered]{algorithm2e}

\title{\LARGE \bf
	A Newton-like Method based on Model Reduction Techniques\\for Implicit Numerical Methods
}

\author{Tobias K. S. Ritschel
\thanks{$^{1}$Tobias K. S. Ritschel is with Department of Applied Mathematics and Computer Science, Technical University of Denmark, DK-2800 Kgs. Lyngby, Denmark
        {\tt\small tobk@dtu.dk}}%
}

\newcommand{\incr}{\,\mathrm{d}}
\newcommand{\pdiff}[2]{\frac{\partial#1}{\partial#2}}

\newtheorem{remark}{Remark}

\definecolor{mygreen}{rgb}{0,0.6,0}
\definecolor{mygray}{rgb}{0.1,0.5,0.5}
\definecolor{mymauve}{rgb}{0.58,0,0.82}
\definecolor{light-gray}{gray}{0.95}

\colorlet{AlgCaptionColor}{light-gray}
\SetKwRepeat{Do}{do}{while}
\makeatletter
\renewcommand{\algocf@makecaption@ruled}[2]{%
	\global\sbox\algocf@capbox{\colorbox{AlgCaptionColor}{\hskip\AlCapHSkip
			\parbox[t]{1.06\hsize}{\algocf@captiontext{\strut#1}{\strut#2\strut}}\hskip.1\algomargin}}
}%
\makeatother
\SetKwInput{KwData}{Input} 
\SetKwInput{KwResult}{Output} 

\begin{document}
	\maketitle
\thispagestyle{empty}
\pagestyle{empty}

	\begin{abstract}
	In this paper, we present a Newton-like method based on model reduction techniques, which can be used in implicit numerical methods for approximating the solution to ordinary differential equations. In each iteration, the Newton-like method solves a reduced order linear system in order to compute the Newton step. This reduced system is derived using a projection matrix, obtained using proper orthogonal decomposition, which is updated in each time step of the numerical method. We demonstrate that the method can be used together with Euler's implicit method to simulate CO$_2$ injection into an oil reservoir, and we compare with using Newton's method. The Newton-like method achieves a speedup of between 39\% and 84\% for systems with between 4,800 and 52,800 state variables.
\end{abstract}

	\section{Introduction}
Newton's method iteratively approximates the solution to a set of nonlinear algebraic equations using successive linearization. It is ubiquitous in numerical methods for simulation, state and parameter estimation, and optimal control of ordinary, partial, differential-algebraic, and stochastic systems of differential equations. In particular, \emph{implicit} numerical methods are suitable for stiff systems, and they require the solution of one or more sets of nonlinear equations in each time step.

Each iteration of Newton's method requires the solution of a linear system of equations. For large-scale systems with many state variables, this is computationally expensive.
Therefore, many variations of Newton-type methods, where the Newton step is approximated, have been proposed~\cite{Morini:1999, Deuflhard:2011}.
Two widely used variations are 1)~the \emph{simplified} Newton method where the Jacobian matrix is not reevaluated in each iteration (such that its factorization can be reused) and 2)~\emph{inexact} Newton methods~\cite{Dembo:etal:1982, Bellavia:etal:2011} where the linear system in each Newton iteration is solved approximately using an iterative method~\cite{Trefethen:Bau:1997, Saad:2003}.

Large-scale systems often arise as the result of spatial discretization of partial differential equations (PDEs)~\cite{Eriksson:etal:1996, LeVeque:2007} or as a network of interacting subsystems~\cite{Porter:Gleeson:2016}. In the former case, it is common that the resulting set of ordinary differential equations (ODEs) is stiff due to fast local changes (compared to the simulation horizon). The stiffness of network systems depends on the dynamics of the individual subsystems and their interactions.
\emph{Model order reduction}, or simply model reduction, methods~\cite{Antoulas:2005, Benner:etal:2017} are relevant to any type of large-scale system. They identify a lower order dynamical system whose state variables can be used to approximate the state of the original dynamical system.
For linear systems, there exist several model reduction methods, and the theory is well-developed.
However, the reduction of general nonlinear systems is an active field of research and remains challenging.

In this work, we propose a Newton-like method where a reduced order linear system is solved in each iteration. We use a projection matrix to compute the system matrices in this reduced system and to compute the approximate Newton step from the solution to the reduced system. We use proper orthogonal decomposition (POD)~\cite{Antoulas:2005} to compute the projection matrix. The method is relevant to implicit numerical methods, and we demonstrate its utility using Euler's implicit method.
A similar method was proposed by Nigro et al.~\cite{Nigro:etal:2016}. However, we compute a projection matrix in each time step whereas they only compute it a few times during the simulation. Instead, they update the projection matrix adaptively during both the time steps and Newton iterations. Finally, if convergence is too slow, we solve the original (full order) linear system for the Newton step. The advantage of the method proposed in this work is its simplicity.
We test the Newton-like method on numerical simulation of CO$_2$ injection into an oil reservoir. This process is modeled using four coupled PDEs which are discretized using a finite volume method.

The remainder of this paper is structured as follows. In Section~\ref{sec:simulation}, we discuss numerical simulation of nonlinear systems using Euler's implicit method and Newton's method, and we discuss model reduction and POD in Section~\ref{sec:model:order:reduction}. In Section~\ref{sec:newton:like:method}, we present the Newton-like method proposed in this paper, and we discuss details of the implementation in Section~\ref{sec:implementation}. Finally, we present the numerical examples in Section~\ref{sec:example}, and we present conclusions in Section~\ref{sec:conclusions}.
	\section{Numerical simulation}\label{sec:simulation}
We consider nonlinear systems in the form
\begin{align}\label{eq:system}
	\dot x(t) &= f(x(t), u(t), d(t), p),
\end{align}
where $x$ are the states, $u$ are manipulated inputs, $d$ are disturbance variables, $p$ are parameters, and $f$ is the right-hand side function. We include the dependency of $f$ on $u$, $d$, and $p$ for completeness, and we assume a zero-order hold parametrization of $u$ and $d$:
\begin{subequations}
	\begin{align}
		u(t) &= u_k, & t \in [t_k, t_{k+1}[, \\
		d(t) &= d_k, & t \in [t_k, t_{k+1}[.
	\end{align}
\end{subequations}

\subsection{Discretization}
We discretize~\eqref{eq:system} using Euler's implicit method:
\begin{align}\label{eq:implicit:euler}
	x_{k+1} - x_k &= f(x_{k+1}, u_k, d_k, p) \Delta t_k.
\end{align}
Here, $x_k \approx x(t_k)$ and $\Delta t_k$ is the $k$'th time step size.
As the right-hand side function, $f$, is nonlinear in the states,~\eqref{eq:implicit:euler} is a set of nonlinear algebraic equations which must be solved for $x_{k+1}$ given $x_k$, $u_k$, $d_k$, $p$, and $\Delta t_k$.

\subsection{Newton's method}
We write the nonlinear algebraic equations~\eqref{eq:implicit:euler} in residual form:
\begin{align}\label{eq:implicit:euler:residual}
	R_k(x_{k+1})
	&= R_k(x_{k+1}; x_k, u_k, d_k, p) \nonumber \\
	&= x_{k+1} - x_k - f(x_{k+1}, u_k, d_k, p) \Delta t_k = 0.
\end{align}
In Newton's method, an approximation of the solution to~\eqref{eq:implicit:euler:residual} is iteratively improved using the update formula
\begin{align}\label{eq:newton:update}
	x_{k+1}^{(\ell+1)} &= x_{k+1}^{(\ell)} + \Delta x_{k+1}^{(\ell)},
\end{align}
where $x_{k+1}^{(\ell)}$ is the $\ell$'th approximation. If
\begin{align}
	\|R_k(x_{k+1}^{(\ell)})\| < \tau
\end{align}
for some tolerance, $\tau$, the iterations are terminated.
The Newton step $\Delta x_{k+1}^{(\ell)}$ is the solution to the linear system of equations obtained by linearizing~\eqref{eq:implicit:euler:residual} around $x_{k+1}^{(\ell)}$:
\begin{align}\label{eq:newton:step:intermediate}
	R_k(x_{k+1}^{(\ell)}) + \pdiff{R_k}{x_{k+1}}(x_{k+1}^{(\ell)}) \Delta x_{k+1}^{(\ell)} &= 0.
\end{align}
The Jacobian matrix is
\begin{align}
	\pdiff{R_k}{x_{k+1}}(x_{k+1}) &= \mathrm I - \pdiff{f}{x}(x_{k+1}, u_k, d_k, p) \Delta t_k,
\end{align}
where $\mathrm I$ is the identity matrix and $\pdiff{f}{x}$ is the Jacobian of $f$.
We rearrange terms in order to write the linear system in the form
\begin{align}\label{eq:newton:step}
	A_k^{(\ell)} \Delta x_{k+1}^{(\ell)} &= b_k^{(\ell)},
\end{align}
where the system matrix and the right-hand side are
\begin{subequations}
	\begin{align}
		A_k^{(\ell)} &= \pdiff{R_k}{x_{k+1}}(x_{k+1}^{(\ell)}), \\
		b_k^{(\ell)} &= -R_k(x_{k+1}^{(\ell)}).
	\end{align}
\end{subequations}
	\section{Model order reduction}\label{sec:model:order:reduction}
In this section, we discuss the challenges of reducing general nonlinear models, and we present the POD method used in the Newton-like method.
The purpose of model reduction is to approximate the states, $x$, by a smaller number of \emph{reduced} states, $\hat x$, e.g., using a linear relation:
\begin{align}\label{eq:model:reduction:approximation}
	x(t) &\approx V \hat x(t).
\end{align}
Here, $x\in \mathbb R^{n_x}$, $V\in\mathbb R^{n_x \times n_r}$, and $\hat x \in \mathbb R^{n_r}$ where the number of reduced states, $n_r$ is significantly smaller than the number of states, $n_x$. Direct substitution of~\eqref{eq:model:reduction:approximation} into the original system~\eqref{eq:system} gives
\begin{align}
	V \dot{\hat x}(t) &= f(V \hat x(t), u(t), d(t), p).
\end{align}
Next, the differential equations are multiplied by the transpose of $W \in \mathbb R^{n_x\times n_r}$ from the left:
\begin{align}
W^T V \dot{\hat x}(t) &= W^T f(V \hat x(t), u(t), d(t), p).
\end{align}
%
%
Assuming that $W^T V$ is invertible, the reduced system is
\begin{align}\label{eq:reduced:system}
	\dot{\hat x}(t) &= \hat f(\hat x(t), u(t), d(t), p),
\end{align}
where the reduced right-hand side function, $\hat f$, is
\begin{align}~\label{eq:reduced:right:hand:side:function}
	\hat f(\hat x(t), &\,u(t), d(t), p) \nonumber \\
	&= (W^T V)^{-1} W^T f(V \hat x(t), u(t), d(t), p).
\end{align}

For general nonlinear systems, evaluating the reduced right-hand side function in~\eqref{eq:reduced:right:hand:side:function} is more expensive than evaluating the right-hand side function in the original system (for linear systems, the reduced system matrices can be computed prior to simulation or analysis). Consequently, using explicit numerical methods for the reduced system~\eqref{eq:reduced:system} is, in most cases, not faster than for the original system. However, for implicit methods, the most computationally expensive step is to solve the linear system of equations in each iteration of the Newton step. Therefore, it is possible to achieve significant speedup by simulating the reduced system. However, it still remains challenging to accurately approximate general nonlinear systems using a reduced system.
\begin{remark}
	The matrix $W^T V$ is often the identity matrix or a diagonal matrix (e.g., for clustering approaches).
\end{remark}
%
%

\subsection{Proper orthogonal decomposition}\label{sec:model:reduction:pod}
In the Newton-like method presented in Section~\ref{sec:newton:like:method}, we use POD~\cite[Sec.~9.1]{Antoulas:2005} to compute the projection matrices $V$ and $W$ based on a matrix of \emph{snapshots}, $X$. That is, each column of $X$ contains the state vector $x_k = x(t_k)$ for some value of $k$. The singular value decomposition of $X$ is
\begin{align}\label{eq:svd}
	X = U \Sigma R,
\end{align}
where $U$ and $R$ are matrices of left and right singular vectors, respectively, and $\Sigma$ is a diagonal matrix with the singular values on the diagonal in descending order. The projection matrices $V = W$ consist of the columns of $U$ corresponding to singular values larger than $\epsilon \Sigma_{11}$. Based on numerical experiments, we choose $\epsilon = 50 \epsilon_m$ where $\epsilon_m = 2^{-52}$ is the machine precision
	\section{Numerical simulation\\using the Newton-like method}\label{sec:newton:like:method}
In the Newton-like method, we replace the Newton update~\eqref{eq:newton:update} by
%
\begin{align}\label{eq:newton:like:update}
	x_{k+1}^{(\ell+1)} &= x_{k+1}^{(\ell)} + V_k \Delta \hat x_{k+1}^{(\ell)},
\end{align}
where $V_k$ is the projection matrix computed using POD, as described in Section~\ref{sec:model:reduction:pod}, and $\Delta \hat x_{k+1}^{(\ell)}$ is the reduced Newton step obtained by solving
\begin{align}\label{eq:newton:like:step}
	\hat A_k^{(\ell)} \Delta \hat x_{k+1}^{(\ell)} &= \hat b_k^{(\ell)}.
\end{align}
We derive the system matrix and the right-hand side following the same steps as in Section~\ref{sec:model:order:reduction}. First, we approximate the Newton step by
\begin{align}
	\Delta x_{k+1}^{(\ell)} &\approx V_k \Delta \hat x_{k+1}^{(\ell)}.
\end{align}
Next, as the linear system would otherwise be overdetermined, we multiply the linear system of equations by $W_k^T$ from the left. The resulting set of linear equations is
%
%
\begin{align}
	W_k^T A_k^{(\ell)} V_k \Delta \hat x_{k+1}^{(\ell)} &= W_k^T b_k^{(\ell)}.
\end{align}
Consequently, the reduced system matrix and right-hand side in~\eqref{eq:newton:like:step} are
\begin{subequations}
	\begin{align}
		\hat A_k^{(\ell)} &= W_k^T A_k^{(\ell)} V_k, \\
		\hat b_k^{(\ell)} &= W_k^T b_k^{(\ell)}.
	\end{align}
\end{subequations}
The dimension of $\hat A_k^{(\ell)}$ is $n_r\times n_r$ which is significantly smaller than the dimension of $A_k^{(\ell)}$ which is $n_x\times n_x$. Therefore, the solution of the linear system~\eqref{eq:newton:like:step} is much less computationally expensive.

\subsection{Algorithms}
We implement Euler's implicit method as shown in Algorithm~\ref{algo:simulation:newton:like:method}. We use Newton's method in the first $N_b$ time steps. For each subsequent time step, we create a snapshot matrix, $X_k$, consisting of up to $N_h$ of the previous time steps. Based on the snapshot matrix, we compute the projection matrices $V_k$ and $W_k$ as described in Section~\ref{sec:model:reduction:pod}. Finally, we use the projection matrices in the Newton-like method to solve the residual equations~\eqref{eq:implicit:euler:residual} for the next time step, $x_{k+1}$. Based on numerical experiments, we choose $N_b = \ln(n_x)$ and $N_h = \sqrt[3]{n_x}$.

We implement the Newton-like method as shown in Algorithm~\ref{algo:newton:like:method}. The initial guess is $x_k$ and the result is $x_{k+1}$. In each iteration, the reduced system~\eqref{eq:newton:like:step} is solved unless the previous reduced Newton step was too small. In line~\ref{line:newton:like:method:remark:1}, we ensure that the full order system is not solved twice in a row. The iterations continue until the norm of the residual equations is below the tolerance $\tau$. For convenience, we also use $\tau$ in line~\ref{line:newton:like:method:remark:2} to determine when to solve the full order system.
\begin{algorithm}[t]
	\KwData{$x_0$, $\{\Delta t_k\}_{k=0}^{N-1}$, $\{u_k\}_{k=0}^{N-1}$, $\{d_k\}_{k=0}^{N-1}$, $p$, \qquad $N$, $N_b$, $N_h$, $\tau$}
	\KwResult{$\{x_k\}_{k=0}^N$}
	\For{$k = 0, \ldots, N-1$}{
		\eIf{$k \geq N_b$}{
			Set $k_h = \max\{0, k - N_h + 1\}$\;
			Set $X_k = \begin{bmatrix} x_{k_h} & \cdots & x_k \end{bmatrix}$\;
			Compute $V_k$ and $W_k$ using POD of $X_k$\;
			Compute $x_{k+1}$ using Algorithm~\ref{algo:newton:like:method}\;
		}{
			Compute $x_{k+1}$ using Newton's method\;
		}
	}
	\caption{Simulation based on Newton-like method}
	\label{algo:simulation:newton:like:method}
\end{algorithm}

\begin{algorithm}[t]
	\KwData{$x_k$, $V_k$, $W_k$, $\tau$}
	\KwResult{$x_{k+1}$}
	Set $\ell = 0$\;
	Set $x_{k+1}^{(0)} = x_k$\;
	\While{$\|R_k(x_{k+1}^{(\ell)})\| \geq \tau$}{
		\eIf{$\ell > 0$ and $\|\Delta \hat x_{k+1}^{(\ell-1)}\| < \tau$}{ \label{line:newton:like:method:remark:2}
			Solve the full order system~\eqref{eq:newton:step} for $\Delta x_k^{(\ell)}$\;
			Compute $x_{k+1}^{(\ell+1)}$ using~\eqref{eq:newton:update}\;
			Set $\|\Delta \hat x_{k+1}^{(\ell)}\|$ larger than $\tau$\; \label{line:newton:like:method:remark:1}
		}{
			Solve the reduced system~\eqref{eq:newton:like:step} for $\Delta \hat x_k^{(\ell)}$\;
			Compute $x_{k+1}^{(\ell+1)}$ using~\eqref{eq:newton:like:update}\;
		}
		Increment $\ell$ by 1\;
	} 
	Set $x_{k+1} = x_{k+1}^{(\ell)}$\;
	\caption{Newton-like method}
	\label{algo:newton:like:method}
\end{algorithm}

	\section{Implementation details}\label{sec:implementation}
We implement Algorithm~\ref{algo:simulation:newton:like:method} and~\ref{algo:newton:like:method} in Matlab~\cite{MATLAB:2021}. The Jacobian matrices of $f$ and $R_k$ are represented as sparse. All other quantities are represented as dense. We use Matlab's \texttt{svds} to compute the left singular vectors corresponding to the $N_h$ largest singular values, and we use Matlab's \texttt{mldivide} (i.e., the backslash operator) to solve the sparse linear systems.
Furthermore, we implement the evaluation of the right-hand side function and its Jacobian using the C++ library DUNE~\cite{Sander:2020} and the thermodynamic library, ThermoLib~\cite{Ritschel:etal:2016, Ritschel:etal:2017b}, and we use a Matlab MEX interface~\cite{Woodford:2022} to call the implementation from Matlab. The MEX interface allocates memory and evaluates the Jacobian of $f$ as sparse before returning it to Matlab, i.e., we do not use Matlab's \texttt{sparse} to convert it from a dense matrix.
Finally, we carry out the computations using Windows Subsystem for Linux 2 (WSL2) on a Windows 10 laptop with 12~MB shared level 3 cache, 1,280~KB level 2 cache for each core, and 48~KB and 32~KB level 1 instruction and data cache for each core, respectively. Furthermore, each core contains two 11'th generation Intel Core i7 3~GHz processors.
	\section{Numerical example}\label{sec:example}
\begin{figure}[t]
	\centering
	\includegraphics[width=0.995\linewidth]{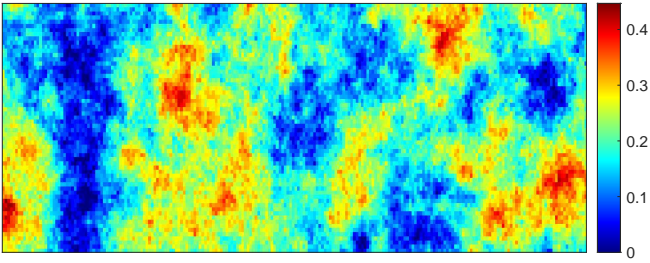}
	\caption{Porosity field.}
	\label{fig:ECC2023_PorosityField}
\end{figure}

\begin{figure}[t]
	\centering
	\includegraphics[width=\linewidth]{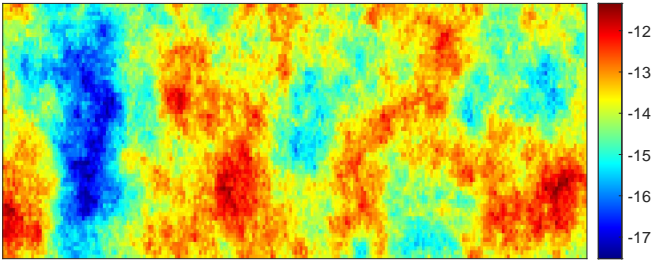}
	\caption{Permeability field in m$^2$.}
	\label{fig:ECC2023_PermeabilityField}
\end{figure}

In this section, we present a numerical example of CO$_2$ injection into an oil reservoir. We provide a brief description of the model in Appendix~\ref{sec:flow:model}, and we use the porosity and permeability fields of the top layer of model 2 from~\cite{Christie:Blunt:2001}. They are shown in Fig.~\ref{fig:ECC2023_PorosityField} and~\ref{fig:ECC2023_PermeabilityField}. Clearly, they are very heterogeneous.
Initially, the reservoir contains water and oil consisting of methane, n-decane and CO$_2$, and the water and oil phases do not mix. Furthermore, we assume that there is no gas phase. We inject liquid CO$_2$, which also contains 1\textperthousand{} methane and 1\textperthousand{} n-decane (mole fractions). The CO$_2$ is injected through two wells located at opposite corners. Each well injects 0.11~m$^3$ per day.

\begin{figure*}[t]
\centering
\includegraphics[width=\textwidth]{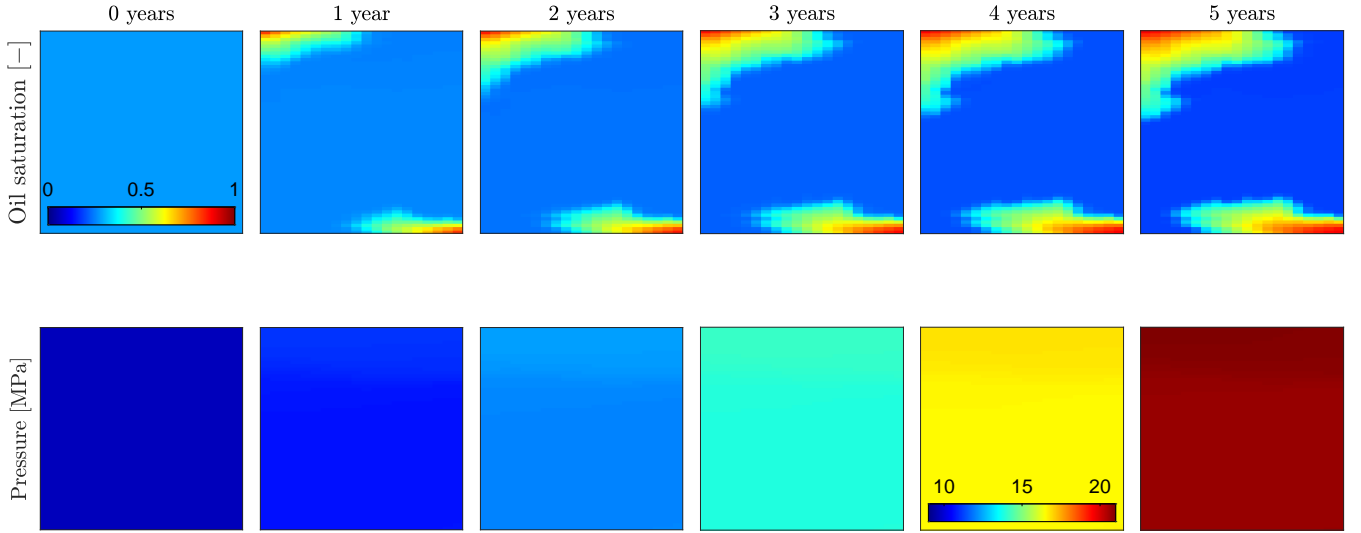}
\caption{Simulation of CO$_2$ injection over 5 years (1,200 grid cells). Top row: Oil saturation. Bottom row: Pressure in MPa.}
\label{fig:ECC2023_Comp_P_6}
\end{figure*}

Fig.~\ref{fig:ECC2023_Comp_P_6} shows an example of a simulation of 1,200 grid cells. The reservoir consists of the left-most 20 columns of cells in the fields shown in Fig.~\ref{fig:ECC2023_PorosityField} and~\ref{fig:ECC2023_PermeabilityField}.
The top row of figures show the oil saturation (ratio between oil and water volume) over time and the bottom row shows the pressure. As CO$_2$ is injected, it displaces water. Therefore, the oil saturation decreases in the middle area. In this simulation, the oil and water cannot leave the reservoir. Consequently, the pressure increases significantly over. Furthermore, the pressure is almost completely uniform in this simulation.

\begin{figure}[t]
	\centering
	\includegraphics[width=\linewidth]{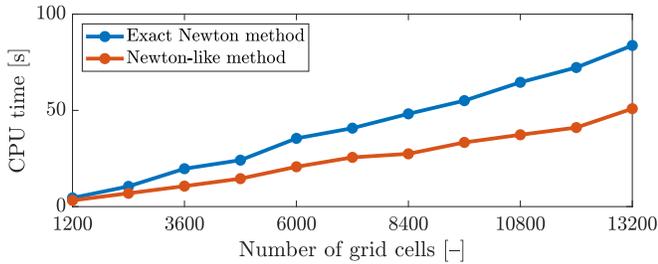}
	\caption{Computation time of simulation using Newton's method and the Newton-like method for different numbers of grid cells.}
	\label{fig:CompareOriginalAndMOR}
\end{figure}

\begin{figure}[t]
	\centering
	\includegraphics[width=0.96\linewidth]{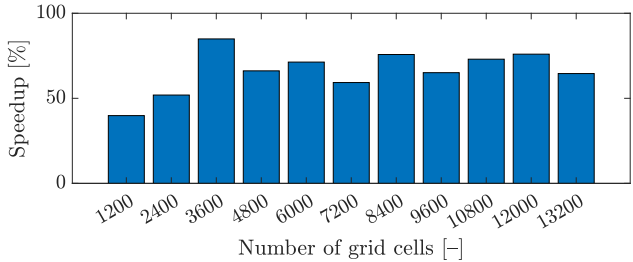}
	\caption{Speedup for the Newton-like method for different numbers of grid cells.}
	\label{fig:CompareOriginalAndMOR_Speedup}
\end{figure}

Fig.~\ref{fig:CompareOriginalAndMOR} shows the computation time of simulating 8~days of CO$_2$ injection using Euler's implicit method and 1)~Newton's method and 2)~the Newton-like method described in Section~\ref{sec:newton:like:method}. We use 40~time steps of 0.1~day and 20~time steps of 0.2~day, and we repeat the simulation for different numbers of grid cells. The number of grid cells changes the simulation scenario because we include a larger part of the reservoir, i.e., we do not refine the discretization. The computation time increases linearly with the number of grid cells, and the increase is lower for the Newton-type method. Fig.~\ref{fig:CompareOriginalAndMOR_Speedup} shows the speedup which is between 39\% (for the lowest number of grid cells) and 84\%. The mean speedup is 66\% and the standard deviation is 12\%.

	\section{Conclusions}\label{sec:conclusions}
In this paper, we propose a Newton-like method based on model reduction methods (specifically, POD). The method can improve the computational efficiency of implicit numerical methods for approximating the solution to general nonlinear ODEs. Such systems are difficult to reduce directly using existing model reduction methods.
The method approximates the linear system solved in each Newton iteration by a lower order system which is significantly less computationally expensive to solve. The system matrix and right-hand side of the reduced linear system are computed based on a projection matrix, which is updated in each time step (except for the first few time steps).
We demonstrate the utility of the method for numerical simulation of CO$_2$ injection into an oil reservoir. We model
this process by four coupled PDEs, and we discretize them using a finite volume method.
Compared to Newton's method, the Newton-like method achieves a speedup of between 39\% and 84\% for 1,200 to 13,200 grid cells, which correspond to 4,800 and 52,800 state variables.

	\appendices

	\section{Compositional flow in porous media}\label{sec:flow:model}
In this appendix, we modify an isothermal and compositional model presented in previous work~\cite{Ritschel:Jorgensen:2019}. Specifically, we assume that there is no gas phase and that the rock and the water are incompressible. The purpose of these simplifications is to obtain a model which (after spatial discretization) consists of a set of ODEs.

The concentrations of water, $C_w$, and the $k$'th component in the oil phase, $C_k$, are described by the PDEs
\begin{subequations}\label{eq:flow:model}
	\begin{align}
		\partial_t C_w &= -\nabla\cdot \mathbf N^w, \\
		\partial_t C_k &= -\nabla\cdot \mathbf N_k + Q_k,
	\end{align}
\end{subequations}
where $\mathbf N^w$ and $\mathbf N_k$ are the molar fluxes of the water component and of the $k$'th component, and $Q_k$ is a source term representing the injection of component $k$.
The molar flux of the $k$'th component is the product of the mole fraction, $x_k$, and the molar flux of the entire oil phase:
\begin{align}
	\mathbf N_k &= x_k \mathbf N^o.
\end{align}
The molar fluxes of both the oil and the water phase are the products of the molar density, $\rho^\alpha$, and the volumetric flux, i.e.,
\begin{align}
	\mathbf N^\alpha &= \rho^\alpha \mathbf u^\alpha, & \alpha &\in \{w, o\},
\end{align}
and the volumetric flux is given by a generalization of Darcy's law to multiphase fluids:
\begin{align}
	\mathbf u^\alpha &= -\frac{k_r^\alpha}{\mu^\alpha} \mathbf K (\nabla P - \bar \rho^\alpha g \nabla z), & \alpha &\in \{w, o\}.
\end{align}
We use Corey's model~\cite{Chen:etal:2006} of the relative permeabilities, $k_r^\alpha$, and we use the model by~Lohrenz et al.~\cite{Lohrenz:etal:1964} to describe the viscosity of the oil phase, $\mu^o$. The viscosity of the water phase is $\mu^w = 0.3$~cP. In this work, we assume the permeability tensor, $\mathbf K$, to be a multiple of the identity matrix, i.e., the permeability is the same in all directions. Furthermore, $P$ is pressure, $\bar \rho^\alpha$ is the mass density, $g$ is the gravity acceleration constant, and $z$ is depth.

As for the flux of component $k$, the source term is the product of the mole fraction and the source term for the entire oil phase:
\begin{align}
	Q_k &= x_k Q^o.
\end{align}
We specify both $x_k$ and $Q^o$.

\subsection{Volume balance and cubic equation of state}
When both oil and gas is present, the phase equilibrium problem for each grid cell is an equality-constrained minimization problem. The objective function is the combined Helmholtz' energy of the rock, water, oil and gas phase and the constraints specify the amounts of moles of each component (in both phases) and the combined volume of all four phases. However, when there is no gas phase, the phase equilibrium problem simplifies to finding a pressure which satisfies the volume balance:
\begin{align}\label{eq:volume:balance}
	V^w + V^o + V^r &= V.
\end{align}
The volumes of the rock, $V^r$, the water, $V^w$, and the entire grid cell, $V$, are independent of the pressure. Consequently, given the amount of moles in the water phase, $n_w$, the oil volume, $V^o$, can be isolated in~\eqref{eq:volume:balance}.
Given the oil volume $V^o = V - V^w - V^r$, we can compute the corresponding molar volume $v^o = V^o/N^o$, where $N^o$ is the total amount of moles in the oil phase, and evaluate the pressure using the cubic equation of state
\begin{align}
	P &= \frac{RT}{v^o - b_m} - \frac{a_m}{(v^o + \epsilon b_m)(v^o + \sigma b_m)}.
\end{align}
Here, $R$ is the gas constant, $T$ is temperature (60$^\circ$C in this work), $a_m = a_m(T, n^o)$ and $b_m(n^o)$ are mixing parameters ($n^o$ is a vector of moles of each component in the oil phase), and $\epsilon$ and $\sigma$ are parameters which depend on the specific cubic equation of state used. We use the Peng-Robinson equation of state.

\begin{remark}
	After spatial discretization, the flow models described previously~\cite{Ritschel:Jorgensen:2019} consist of differential-algebraic equations (DAEs) in the form
	\begin{subequations}
		\begin{align}
			\dot x(t) &= F(y(t), u(t), d(t), p), \\
			0 &= G(x(t), y(t), z(t), p),
		\end{align}
	\end{subequations}
	whereas the model described above is in the form
	\begin{subequations}
		\begin{align}
			\dot x(t) &= F(y(t), u(t), d(t), p), \\
			\label{eq:flow:model:explicit:algebraic:equations}
			y(t) &= \tilde G(x(t), p).
		\end{align}
	\end{subequations}
	Here, $y$ and $z$ are algebraic variables, $F$ is the right-hand side function, and $G$ and $\tilde G$ are algebraic functions. As the algebraic equations~\eqref{eq:flow:model:explicit:algebraic:equations} are explicit in the algebraic variables, the above model is also in the form~\eqref{eq:system}, i.e.,
	\begin{align}
		\dot x(t)
		&= f(x(t), u(t), d(t), p) \nonumber \\
		&= F(\tilde G(x(t), p), u(t), d(t), p).
	\end{align}
\end{remark}

\subsection{Discretization}
The PDEs~\eqref{eq:flow:model} are in the form
\begin{align}
	\partial_t C &= -\nabla\cdot \mathbf N + Q,
\end{align}
and we discretize them using a finite volume method. First, we integrate over each grid cell, $\Omega_i$:
\begin{align}
	\partial_t \int_{\Omega_i} C \incr V &= -\int_{\Omega_i} \nabla\cdot \mathbf N \incr V + \int_{\Omega_i} Q \incr V.
\end{align}
The left-hand side can be evaluated exactly, i.e.,
\begin{align}
	\int_{\Omega_i} C \incr V &= n_i,
\end{align}
where $n_i$ is the total amount of moles in the $i$'th grid cell. Next, we use Gauss' divergence theorem, write the surface integral as the sum of the integrals over each face of the grid cell, $\gamma_{ij}$, and approximate the integrals:
\begin{align}\label{eq:flux:gauss:divergence}
	\int_{\Omega_i} \nabla\cdot \mathbf N \incr V
	&= \int_{\partial \Omega_i} \mathbf N \cdot \mathbf n \incr A \nonumber \\
	&= \sum_{j\in\mathcal N^{(i)}} \int_{\gamma_{ij}} \mathbf N \cdot \mathbf n \incr A \nonumber \\
	&\approx \sum_{j\in\mathcal N^{(i)}} (A \mathbf N \cdot \mathbf n)_{ij}.
\end{align}
The set $\mathcal N^{(i)}$ contains the indices of the grid cells that are adjacent to the $i$'th grid cell, $A_{ij}$ is the area of face $\gamma_{ij}$ (the face shared by $\Omega_i$ and $\Omega_j$), and $\mathbf n$ is the outward normal vector on the boundary of the grid cell, $\partial \Omega_i$.
Similarly, we approximate the integral of the source term as
\begin{align}
	\int_{\Omega_i} Q \incr V & \approx (Q V)_i.
\end{align}
The resulting approximation in~\eqref{eq:flux:gauss:divergence} contains the flux evaluated at the center of each face of the grid cell. We use a two-point flux approximation~\cite{Lie:2019} to approximate this flux. The result is
\begin{subequations}
	\begin{align}
		(A \mathbf N^w \cdot \mathbf n)_{ij} &\approx -(\Gamma \hat H^w \Delta \Phi^w)_{ij}, \\
		(A \mathbf N_k \cdot \mathbf n)_{ij} &\approx -(x_k \Gamma \hat H^o \Delta \Phi^o)_{ij},
	\end{align}
\end{subequations}
where $\Gamma_{ij}$ is the geometric part of the transmissibilities, $\hat H^\alpha$ is the fluid part of the transmissibilities, and $\Delta \Phi^\alpha$ is difference in potential. The geometric part of the transmissibilities is given by
\begin{subequations}
	\begin{align}
		\Gamma_{ij} &= A_{ij} \left(\hat \Gamma_{ij}^{-1} + \hat \Gamma_{ji}^{-1}\right)^{-1}, \\
		\hat \Gamma_{ij} &= \left(\mathbf K_i \frac{c_{ij} - c_i}{|c_{ij} - c_i|^2}\right) \cdot \mathbf n,
	\end{align}
\end{subequations}
where $c_{ij}$ is the center of $\gamma_{ij}$ and $c_i$ is the center of $\Omega_i$. Furthermore, the difference in potential is
\begin{align}
	\Delta \Phi_{ij}^\alpha &= (\Delta P - \bar \rho^\alpha g \Delta z)_{ij},
\end{align}
where we approximate the mass density at the face center using an average, i.e.,
\begin{align}
	\bar \rho_{ij}^\alpha &= (\bar \rho_i^\alpha + \bar \rho_j^\alpha)/2,
\end{align}
and the differences in pressure and depth are
%
	\begin{align}
		\Delta P_{ij} &= P_j - P_i, &
		\Delta z_{ij} &= z_j - z_i.
	\end{align}
%
Finally, the fluid part of the transmissibilities is upwinded based on the difference in potential:
\begin{align}
	\hat H_{ij}^\alpha &=
	\begin{cases}
		(\rho^\alpha k_r^\alpha/\mu^\alpha)_i & \text{if}~\Delta \Phi_{ij}^\alpha <    0, \\
		(\rho^\alpha k_r^\alpha/\mu^\alpha)_j & \text{if}~\Delta \Phi_{ij}^\alpha \geq 0.
	\end{cases}
\end{align}

	\bibliography{./ref/ref.bib}

\begin{thebibliography}{10}

\bibitem{Morini:1999}
B.~Morini, ``Convergence behaviour of inexact {N}ewton methods,'' {\em
  Mathematics of Computation}, vol.~68, no.~228, pp.~1605--1613, 1999.

\bibitem{Deuflhard:2011}
P.~Deuflhard, {\em Newton methods for nonlinear problems: Affine invariance and
  adaptive algorithms}.
\newblock Springer Series in Computational Mathematics, Springer, 2011.

\bibitem{Dembo:etal:1982}
R.~S. Dembo, S.~C. Eisenstat, and T.~Steihaug, ``Inexact {N}ewton methods,''
  {\em SIAM Journal on Numerical Analysis}, vol.~19, no.~2, pp.~400--408, 1982.

\bibitem{Bellavia:etal:2011}
S.~Bellavia, S.~Magheri, and C.~Miani, ``Inexact {N}ewton methods for model
  simulation,'' {\em International Journal of Computer Mathematics}, vol.~88,
  no.~14, pp.~2969--2987, 2011.

\bibitem{Trefethen:Bau:1997}
L.~N. Trefethen and D.~Bau, {\em Numerical Linear Algebra}.
\newblock SIAM, 1997.

\bibitem{Saad:2003}
Y.~Saad, {\em Iterative methods for sparse linear systems}, vol.~2.
\newblock SIAM, 2003.

\bibitem{Eriksson:etal:1996}
K.~Eriksson, D.~Estep, P.~Hansbo, and C.~Johnson, {\em Computational
  differential equations}.
\newblock Studentlitteratur, 1996.

\bibitem{LeVeque:2007}
R.~J. LeVeque, {\em Finite difference methods for ordinary and partial
  differential equations: Steady-state and time-dependent problems}.
\newblock SIAM, 2007.

\bibitem{Porter:Gleeson:2016}
M.~A. Porter and J.~P. Gleeson, {\em Dynamical systems on networks: A
  tutorial}, vol.~4 of {\em Frontiers in Applied Dynamical Systems: Reviews and
  Tutorials}.
\newblock Springer, 2016.

\bibitem{Antoulas:2005}
A.~C. Antoulas, {\em Approximation of large-scale dynamical systems}.
\newblock Advances in Design and Control, SIAM, 2005.

\bibitem{Benner:etal:2017}
P.~Benner, A.~Cohen, M.~Ohlberger, and K.~Willcox, eds., {\em Model reduction
  and approximation: Theory and algorithms}.
\newblock Computational Science \& Engineering, SIAM, 2017.

\bibitem{Nigro:etal:2016}
P.~S.~B. Nigro, M.~Anndif, Y.~Teixeira, P.~M. Pimenta, and P.~Wriggers, ``An
  adaptive model order reduction with quasi-{N}ewton method for nonlinear
  dynamical problems,'' {\em International Journal for Numerical Methods in
  Engineering}, vol.~106, pp.~740--759, 2016.

\bibitem{MATLAB:2021}
MATLAB, {\em version 9.11.0 (R2021b)}.
\newblock Natick, Massachusetts: The MathWorks Inc., 2021.

\bibitem{Sander:2020}
O.~Sander, {\em {DUNE} -- {T}he distributed and unified numerics environment},
  vol.~140 of {\em Lecture Notes in Computational Science and Engineering}.
\newblock Springer, 2020.

\bibitem{Ritschel:etal:2016}
T.~K.~S. Ritschel, J.~Gaspar, A.~Capolei, and J.~B. J{\o}rgensen, ``An
  open-source thermodynamic software library,'' Tech. Rep. DTU Compute
  Technical Report-2016-12, 2016.

\bibitem{Ritschel:etal:2017b}
T.~K.~S. Ritschel, J.~Gaspar, and J.~B. J{\o}rgensen, ``A thermodynamic library
  for simulation and optimization of dynamic processes,'' {\em
  IFAC-PapersOnLine}, vol.~50, no.~1, pp.~3542--3547, 2017.

\bibitem{Woodford:2022}
O.~Woodford, ``Example {MATLAB} class wrapper for a {C}++ class.''
  www.mathworks.com/matlabcentral/fileexchange/38964-example-matlab-class-wrapper-for-a-c-class.
\newblock MATLAB Central File Exchange. Retrieved November 28, 2022.

\bibitem{Christie:Blunt:2001}
M.~A. Christie and M.~J. Blunt, ``Tenth {SPE} comparative solution project: A
  comparison of upscaling techniques,'' {\em SPE Reservoir Evaluation \&
  Engineering}, vol.~4, no.~4, pp.~308--317, 2001.

\bibitem{Ritschel:Jorgensen:2019}
T.~K.~S. Ritschel and J.~B. J{\o}rgensen, ``Dynamic optimization of
  thermodynamically rigorous models of multiphase flow in porous subsurface oil
  reservoirs,'' {\em Journal of Process Control}, vol.~78, pp.~45--56, 2019.

\bibitem{Chen:etal:2006}
Z.~Chen, G.~Huan, and Y.~Ma, {\em Computational methods for multiphase flow in
  porous media}.
\newblock Computational Science \& Engineering, SIAM, 2006.

\bibitem{Lohrenz:etal:1964}
J.~Lohrenz, B.~G. Bray, and C.~R. Clark, ``Calculating viscosities of reservoir
  fluids from their compositions,'' {\em Journal of Petroleum Technology},
  vol.~16, no.~10, pp.~1171--1176, 1964.

\bibitem{Lie:2019}
K.-A. Lie, {\em An introduction to reservoir simulation using {MATLAB}/{GNU}
  {O}ctave: User guide for the {MATLAB} reservoir simulation toolbox ({MRST})}.
\newblock Cambridge University Press, 2019.

\end{thebibliography}

\end{document}